\documentclass[12pt]{article}
\usepackage[all]{xy}
\usepackage{amsthm,amsfonts}
\xyoption{all}
\newtheorem{teo}{Theorem}[section]
\newtheorem{lema}[teo]{Lemma}
\newtheorem{prop}[teo]{Proposition}

\newtheorem{cor}[teo]{Corollary}

\newcommand{\g}{\overline{\mathcal{H}}_g}
\newcommand{\binf}{\overline{\mathcal{B}}_{2g+2}}

\newcommand{\ohyp}{\mathcal{H}_g}
\newcommand{\bin}{\mathcal{B}_{2g+2}}
\newcommand{\dgd}{\overline{\mathcal{M}}_{0,2g+2}}
\newcommand{\odgd}{\mathcal{M}_{0,2g+2}}

\newcommand{\n}{\overline{\mathcal{M}}_{0,n}}
\newcommand{\bn}{\overline{\textbf{M}}_{0,n}}
\def\p#1{\mbox{{\rm I$\!$P}$_{#1}$}}
\def\P#1{\mathbb{P}_#1}

\def\complex{{\rm C}\hskip-6pt\raise3pt\hbox{$_{^|}$}~}
\def\QED{\ifmmode\squareforqed\else{\unskip\nobreak\hfil
        \penalty50\hskip1em\null\nobreak\hfil\squareforqed
        \parfillskip=0pt\finalhyphendemerits=0\endgraf}\fi
                \medskip}
\begin{document}

 \hoffset = -1truecm
\voffset = -2truecm

\title{\textbf{  The Moduli Spaces of 
 Hyperelliptic Curves and Binary Forms}}
  \author{
\textbf{ D. Avritzer\thanks{Both authors would like
to thank the GMD (Germany) and CNPq (Brasil)
for support during the preparation of this
paper}}
\and \textbf{  H. Lange}$^{*}$}
\date{}
\maketitle

\section{Introduction}

Every smooth hyperelliptic curve of genus $g \ge 2$ can be considered as a 
binary form of degree $2g+2$ with nonzero discriminant in an obvious way.
So if $\ohyp$   and $\bin$ denote the coarse moduli spaces
of smooth hyperelliptic curves of genus $g$
and binary forms of degree $2g+2$ with non vanishing discriminant
respectively, there is a canonical isomorphism 
\[ \ohyp \cong \bin . \]
Let $\g$ and $\binf$ denote the moduli spaces of stable hyperelliptic curves
of genus $g$  and semistable binary forms of degree $2g+2$ respectively.
These are certainly different compactifications of $\ohyp=\bin.$
In fact the boundary $\Delta=\g\setminus \ohyp$ is a divisor whereas
$\binf$ is a one point compactification of the moduli space 
of stable binary forms.
Moreover, $\g$ and $\binf$  are constructed as quotients
via different group actions: $\g$ is defined as the closure of $\ohyp$ in
$\overline\mathcal{M}_g$ the moduli space of stable curves of genus $g$ 
and this is constructed
using the group $PGL(6g-5),$ whereas $\binf$ is constructed classically
using the group $SL_2(\complex).$
It is the aim of this note to work out the relation
between the spaces $\g$ and $\binf.$
The main result is the following theorem, which generalizes Theorem 5.6
of \cite{kn:pqh} where we give a proof in the special case $g=2.$
\bigskip

\begin{teo} \label{princ}
The canonical isomorphism $
\ohyp \widetilde{\longrightarrow } \bin$
extends to a holomorphic map $f_g: \g \longrightarrow \binf .$
\end{teo}

Moreover we work out how the boundary components $\Delta_i$
and $\Xi_i$ of $\Delta$ are contracted under the map $f_g.$
\bigskip

Turning to a more detailed description, we study, in Section 2, the
relation of $\g$ to some other moduli spaces. We show that there are
canonical isomorphisms of $\g$ to the moduli spaces
$\mathbf{\overline{H}_{2,g}}$ of admissible double covers of $(2g+2)-$marked 
curves of genus zero  and to the moduli space
$\overline{\mathcal{M}}_{0,2g+2}$ of stable $(2g+2)-$marked curves
of genus 0. For the definition of these spaces see Section 2.
Certainly the results of Section 2 are well known to the specialists.
In Section 3, we show that for any $ m \geq 3$ there 
there is a canonical holomorphic map
$\overline{\mathcal{M}}_{0,m} \longrightarrow \overline\mathcal{B}_m.$
The main point for this is that any stable $m-$marked
curve $C$ of genus 0 outside a certain irreducible boundary component
of $\overline\mathcal{M}_{0,m}$ admits a central component and for this fact we give a combinatorial
proof using the dual graph of $C.$
In Section 4, we determine the  stable reduction of any curve $y^2=f(x)$
where the homogenization of $f(x)$ is a stable binary form
of degree $2g+2.$
This is applied in Section 5 to describe the holomorphic map
$f_g: \g \longrightarrow \binf$ explicitly  and in particular to work out
how the boundary divisors of $\g$ are contracted.

\section{The Moduli Space of Hyperelliptic Curves}

Let $\mathcal{M}_g$ denote the moduli space of smooth curves of genus 
$g \geq 2$ over the field of complex numbers.
Let $\overline{\mathcal{M}}_g$ denote
 the Deligne-Mumford compactification
 consisting of {\it stable} curves of arithmetic genus 
$g$, that is, curves whose only singularities are nodes and whose rational 
components contain at least three singular points of the curve. 
Curves whose only singularities are nodes and whose rational components
contain at least two singular points of the curve are called {\it semistable.}

\bigskip

Let $\ohyp$ denote the moduli space of smooth hyperelliptic curves
of genus $g,$ considered as a subspace of  $\overline{\mathcal{M}}_g.$
The curves representing elements of the closure $\g$ of $\ohyp$ in 
$\overline{\mathcal{M}}_g$ are called {\it stable 
 hyperelliptic curves.}

In order to study stable hyperelliptic curves it is helpful to introduce
the concept of admissible covers over $n-$marked curves of genus zero. 

For this, recall that
{\it a stable} (respectively {\it semistable) $n-$pointed curve} 
is by definition a complete connected curve $B$ 
that has only nodes as singularities,
together with an \textbf{ ordered} collection $p_1,\dots, p_n \in B$ of 
distinct smooth points such that
every smooth rational component of the normalization of $B$ has at least
3 (respectively 2) points lying over singular points
or points among $p_1,\dots,p_n.$ Note that an $n-$pointed
curve of genus 0 is stable if and only if it admits no nontrivial automorphisms.
Similarly we define
{\it a stable} (respectively {\it semistable) $n-$marked curve}
 to be  a complete connected curve $B$ that has only nodes as singularities,
together with an \textbf{ unordered} collection $p_1,\dots, p_n \in B$ 
of distinct smooth points such that
every smooth rational component of the normalization 
of $B$ has at least 3 (respectively 2)  points lying over singular points
of $B$ or points among $p_1,\dots,p_n.$
Note that any stable $n-$marked curve may admit nontrivial automorphisms,
but its group of automorphisms is always finite. In the genus 0 case
 it is  a subgroup of the
symmetric group of order $n.$

\begin{prop}
For $n \geq 3,$
the  coarse moduli space $\n$ of stable $n-$marked curves of genus 0 exists.
\end{prop}

\textsc{Proof:} According to \cite{kn:kn}, the coarse  moduli space 
$\bn$  of stable $n-$pointed curves of genus 0
exists. It is even a fine moduli space.
The symmetric group $S_n$ of degree $n$  acts on $\bn$
in an obvious way. The group $S_n$ being finite, it is easy to see
that the quotient
$\overline{\mathcal{M}}_{0,n}={\overline\textbf{M}}_{0,n}/S_n$ 
is a coarse moduli space of stable $n-$marked curves of  genus 0. \qed

\bigskip

Let $(B,p_1,\dots,p_n) \in \overline{\mathcal{M}}_{0,n}$ be a stable 
n-marked curve of genus $0$ and let $q_1,\dots,q_k$ 
be the nodes of the curve $B.$ 
By an {\it admissible $d-$fold cover }of the curve $(B,p_1,\dots,p_n)$
we  mean a connected nodal curve $C$ together with a regular map 
$\pi: C \longrightarrow B$
such that:\newline
1. $\pi^{-1}(B_{ns})=C_{ns}$ and the restriction of the map $\pi$ to this open set
is a $d-$fold cover simply branched over the points $p_i$ 
and otherwise unramified; and \newline
2. $\pi^{-1}(B_{sing})=C_{sing}$ and for every node $q$ of $B$ 
and every node $r$ of
$C$ lying over it, the two branches of $C$ near $r$ map to the branches of $B$ 
near $q$ with the same ramification index.

It is clear when two admissible covers should be called isomorphic.

\begin{prop}
For any $g \geq 1$ the coarse moduli space 
{\rm $\overline{\textbf{H}}_{d,g},$} of admissible $d-$fold covers of stable $(2(g+d)-2)-$marked curves of genus 0, exists.
\end{prop}
\bigskip

\textsc{Proof:}
In  \cite{kn:mum} the analogous notion of admissible covers of a stable
$n-$pointed curve of genus zero was defined and it was shown
that the coarse moduli space $\overline{\mathcal{H}}_{d,g}$ of admissible 
 $d-$fold covers of stable $(2(g+d)-2)-$pointed curves of genus 0 exists. 
The symmetric group $S_{2(g+d)-2}$ of $2(g+d)-2$ letters  acts 
in an obvious way
and since it is finite it is easy to see that the quotient 
$\overline{\textbf{H}}_{d,g}=\overline{\mathcal{H}}_{d,g}/S_{2(d+g)-2}$
satisfies the assertion. \qed    
\bigskip

We need the following properties of admissible double covers.

\begin{lema} \label{lema1}
Let $\pi:C \longrightarrow B$ be an admissible  double cover of a curve 
$(B,p_1,\dots,p_{2n}) \in \overline{\mathcal{M}}_{0,2n}.$ \newline
{\rm (a)} Every component of $C$ is smooth. \\
{\rm (b)} $C$ is a semistable curve.\\
{\rm (c)} The stable reduction $C_s$ of $C$ is obtained by contracting
those rational components of $C$ which intersect the other components 
of $C$ in only 2 points.\\
{\rm (d)} If $C_i$ is a rational component of $C$ intersecting the other
components only in $x_1$ and $x_2$ then $\pi(x_1)$=$\pi(x_2).$ 
\end{lema}

\textsc{Proof:} 
(a) Every component of $B$ is smooth, in fact isomorphic to \p{1},
since it is of arithmetic genus 0. This implies (a), 
the map $\pi$ being a regular double cover.\\
(b) Suppose $C$ is not semistable, so $C$ contains a rational component $C_i$
intersecting the other components of $C$ in only one point, $x$ say.
Then $x$ is a branch point of the cover $\pi$ and $B_i=\pi(C_i)$ 
contains only one of the marked points, 
namely the other ramification point.
 This contradicts the stability of $B.$\\
(c) and (d) Suppose $C_i$ is a component of $C$ 
violating the stability of $C.$ Then $C_i$ is a rational component 
intersecting the other components
of $C$ in only 2 points, $x_1$ and $x_2$ say. 
Suppose $\pi(x_1) \neq \pi(x_2)$ and let $B_i=\pi(C_i).$
So $\pi$ is ramified at $x_1$ and $x_2.$
It follows that $\pi|_{C_i}:C_i \longrightarrow B_i$ does not contain
any marked points contradicting the stability of $B.$ \qed

\begin{prop}\label{iso}
There is a canonical isomorphism {\rm $\varphi:\overline{\textbf{H}}_{2,g}
\longrightarrow \overline{\mathcal{H}}_g$}
of the moduli space {\rm $\overline{\textbf{H}}_{2,g}$} of admissible double covers 
of stable $(2g+2)-$marked curves
of genus 0 onto the moduli space $\g$
of stable hyperelliptic curves of genus $g.$
\end{prop}

\textsc{Proof:} Let $\pi : C \longrightarrow B$ be an admissible
double cover in   $\overline{\textbf{H}}_{2,g}.$ It is certainly a limit 
of a family of smooth double covers of \p{1}. Hence $C$ is a limit 
of hyperelliptic curves. Define $\varphi( \pi:C \longrightarrow B)$
to be the stable reduction $C_s$ of $C.$ Since this map can be 
defined for families of admissible double covers , 
we obtain a morphism $\varphi:
\overline{\textbf{H}}_{2,g} \longrightarrow \overline{\mathcal{H}}_g.$

To define an inverse morphism, let $C$ be a stable hyperelliptic curve.
According to (\cite{kn:HM}), Theorem 3.160, there is an admissible double cover
$\pi:C'\longrightarrow B$ of a stable marked curve $B$ of genus $0,$
such that $C'$ is stably equivalent to $C.$
It follows from Lemma ~\ref{lema1} that $\pi: C' \longrightarrow B$ 
is uniquely determined by $C$ and is in fact the blow up of the 
singular points of the components of $C$ as well as all points were 2 curves of genus 0 meet.
In this way we obtain a map $\psi: \overline{\mathcal{H}}_g \longrightarrow
\overline{\textbf{H}}_{2,g}.$
Since this definition extends to families of stable maps, $\psi$ is a morphism.
By construction $\varphi$ and $\psi$ are inverse of each other.\qed

\begin{cor} \label{iso1}
There is a canonical isomorphism
 \[ \Phi :\overline{\mathcal{H}}_g \longrightarrow \overline{\mathcal{M}}_{0,2g+2}\]
of the moduli space $\g$ of stable hyperelliptic curves of genus $g$
onto the moduli space $\dgd$ of stable $(2g+2)-$marked curves of genus
0.
\end{cor}
\textsc{Proof:} 
By Proposition \ref{iso} it suffices to show that there is a canonical 
isomorphism $\varphi: \overline{\textbf{H}}_{2,g} \longrightarrow \dgd .$
If $\pi: C \longrightarrow B$ represents an element of 
$\overline{\textbf{H}}_{2,g},$ the definition of admissible double 
cover and the Hurwitz formula
imply that $B$  is a stable $(2g+2)-$marked curve of arithmetic genus $0.$ 
This gives a map $\varphi: \overline{\textbf{H}}_{2,g} \longrightarrow \dgd$
which certainly is a morphism. Conversely, in order to define the map 
$\psi: \dgd \longrightarrow \overline{\textbf{H}}_{2,g},$ 
let 
$B$ be a stable $(2g+2)-$marked curve of genus 0. We want to show we can 
produce from $B$ in a unique way  a nodal curve $C$ and a regular map 
$\pi$ such that
$\pi:C \longrightarrow B$ is an admissible 2-sheeted cover corresponding to
$B.$ 
Write $B=L_1\cup \dots \cup L_k$ where the $L_i$ are the irreducible 
components of $B.$ Let $p_1,\dots p_{k-1}$ be the nodes of $B$
and $\alpha_1^j,\dots,\alpha_{i_j}^j$
be the marks on $L_j.$ Let $L_1,\dots,L_s$ be the components
of $B$ that intersect the other components in only one point.
Reenumerating the $L_i$ if necessary we may assume $p_1,\dots,p_s$ are
the only nodes of the components 
 $L_1,\dots,L_s.$
We have the following possibilities for $C_j, j=1,\dots, s:$ \newline
1. If ${i_j}$ is even then let $C_j$ be the unique
double cover of $L_j$ branched over the points $\alpha_i^j,j=1,\dots,i_j.$
The curve $C_j$ will  not be ramified over $p_j.$\newline
2. If ${i_j}$ is odd then let $C_j$ be the unique
 double cover of $L_j$ ramified over the $\alpha_i^j,i=1,\dots,i_j$ 
and also over $p_j.$\newline
Now we consider the curve $\overline{B\setminus (L_1\cup\dots\cup L_s)},$
if it is nonempty, 
adding one mark
for each node that is a ramification point in the first step 
and repeat the process
above determining in a unique way the double covers $C_i$ over the components
$L_i$ ramified in the marked points of each component.
The admissible double cover, the image of $B$ under $\psi$ will be the curve $C$
given by $C=C_1\cup \dots\cup C_k.$ It is easy to see that $\psi$ is a morphism.
By construction the maps $\varphi$ and $\psi$ are inverse of each other
and we have the result. 
\qed

\bigskip
The isomorphism $\Phi$ maps the open set $\ohyp$ of smooth hyperelliptic 
curves in $\g$ onto the open set $\odgd$ of smooth $(2g+2)-$marked curves
of genus $0$ in $\dgd$. So the map $\Phi$ can be used to determine
the boundary $\g \setminus \ohyp$ from the boundary  $\dgd \setminus \odgd.$
In fact, Keel shows in \cite{kn:ke} that the boundary  $\dgd \setminus \odgd$
is a divisor consisting of $g+1$ irreducible divisors $D_i$
\[ \dgd \setminus \odgd = D_1 \cup \dots \cup D_{g+1} \]
where a general point of $D_j$ represents a curve $B=B_1 \cup B_2$
with $B_1$ and $B_2$ isomorphic to \p{1} intersecting transversely 
in one point $p$ and $B_1$ containing exactly $j$ points of the marking.
Define
\[ \Delta_i:=\Phi^{-1}(D_{2i+1}) \mbox{ for }i=1,\dots,[ \frac{g}{2} ] \]

and

\[ \Xi_i:=\Phi^{-1}(D_{2i+2}) \mbox{ for } i=0,\dots,[\frac{g-1}{2}]. \]
Here we follow the notation of Cornalba-Harris (see \cite{kn:HC})
who computed the boundary $\g \setminus \ohyp$ in a different way.
Then it is clear from the definition of the map $\Phi$ that a general member
of $\Delta_i$ is represented by an admissible covering as shown in Figure 1
whereas a general element of $\Xi_i$ is represented by an admissible
covering as shown in Figure 2. Note that the curve $C=C_1 \cup C_2$ 
is always stable except in case of $\Xi_0$ where $C_1$ has to be contracted
to yield an irreducible rational curve with a node.

\begin{figure*}
\parbox{3in}{
\begin{picture}(200,200)
{\thicklines
\put(10,10){\line(4,1){110}}
\put(180,10){\line(-4,1){110}}
\put(10,130){\oval(170,30)[r]}
\put(180,130){\oval(170,30)[l]}}
\put(90,155){$q$}
\put(20,160){$C_1$}
\put(160,160){$C_2$}
\put(20,20){$B_1$}
\put(160,20){$B_2$}
\put(20,60){{\tiny 2i+1 branch points}}
\put(110,60){{\tiny 2g-2i+1 branch points}}
\put(90,15){$p$}
\end{picture}
\caption{A general member of $\Delta_i$}}
\parbox{3in}{
\begin{picture}(200,200)
{\thicklines
\put(10,10){\line(4,1){110}}
\put(180,10){\line(-4,1){110}}
\put(5,130){\oval(190,30)[r]}
\put(187,130){\oval(190,30)[l]}}
\put(90,160){$q_1$}
\put(20,160){$C_1$}
\put(160,160){$C_2$}
\put(20,20){$B_1$}
\put(160,20){$B_2$}
\put(20,60){{\tiny 2i+2 branch points}}
\put(120,60){{\tiny 2g-2i branch points}}
\put(90,100){$q_2$}
\put(90,15){$p$}
\end{picture}
\caption{A general member of $\Xi_i$}}
\end{figure*}

\section{Marked curves of genus 0 and binary forms}

A {\it binary form of degree $m$ } is by definition a homogeneous
polynomial $f(x,y)$ of degree $m$ in 2 variables $x$ and $y$ over the field
of complex numbers. We consider binary forms only up to a multiplicative
constant. Hence we consider a binary form  as a smooth "$m-$marked"
curve $(\p{1},x_1,\dots,x_m)$ of genus 0. The "marking"
$x_1,\dots,x_m$ is given by the roots of the form $f$ counted with 
multiplicities. Recall that a binary form of degree $m$ is called 
{\it stable} (respectively {\it semistable)} if 
no  root of $f$ is of multiplicity $\geq \frac{m}{2}$ 
(respectively $> \frac{m}{2}$).
According to (\cite{kn:MF}) the moduli space 
$\overline{\mathcal{B}}_m$ of (equivalence classes of) semistable 
binary forms of degree $m$ exists and is a projective variety
of dimension $m-3.$ Note that for $m$ odd there are no
strictly semistable binary forms. It is well known 
(see (\cite{kn:GEY}) that for $m=2n$ even  all
strictly semistable binary forms of degree $2n$ correspond to one point in 
$\overline{\mathcal{B}}_{2n}$ called {\it the semistable point} 
of $\overline{\mathcal{B}}_{2n}.$

The elements of the open dense set $\mathcal{M}_{0,m}$ of  
$\overline{\mathcal{M}}_{0,m}$ representing a smooth $m-$marked curve may 
be considered as a binary form with no multiple
roots. This induces a morphism

\[ \mathcal{M}_{0,m} \longrightarrow \overline{\mathcal{ B}}_{m}  \]

 It is the aim of this section to prove the following

\begin{teo} \label{main}
The map $\mathcal{M}_{0,m}\longrightarrow \overline{\mathcal{B}}_{m}$ 
extends to a holomorphic map $F_{m}: \overline{\mathcal{M}}_{0,m} 
\longrightarrow \overline{\mathcal{B}}_{m}$ for every $m \geq 3.$
\end{teo}

The idea of the proof is to show that for any stable $m-$marked
curve of genus 0 there is a unique "central component", to which 
one can contract all branches in order to obtain a stable binary form.
For this we need some preliminaries.
\bigskip

Recall that to any nodal curve $C$ one can associate its dual graph:
each vertex of the graph represents a component of $C$ and two vertices
are connected by an edge if the corresponding components intersect.
Thus a connected tree all of whose vertices  represent smooth curves
of genus 0 corresponds to a nodal curve of genus 0. Since we consider
only connected curves, every tree will be connected without further
saying. To an $m-$marked nodal curve $(C,p_1,\dots,p_m)$ of genus 0
we associate a {\it weighted tree} in the following way:
Let $T$ denote the  tree associated to $C$ and $v_i$ a vertex 
representing a component $C_i$ of $C$ which contains exactly $n_i$ 
of the points $p_1,\dots,p_m.$ Then we associate the weight $wt(v_i):=n_i$
to the vertex $v_i.$
The number $m$ of marked  points will be called the {\it total weight}
of the tree $T.$ If $v \in T$ is  a vertex, we denote by $e_j(v)$
the edges with one end in $v$ and call them the edges {\it starting at} $v.$
A (connected) weighted tree is called {\it stable}  if for every 
vertex $v$ the sum of its weight and the number of edges starting
at $v$ is $\geq 3.$
Thus a marked nodal curve of genus 0 is stable if and only if
its corresponding tree $T$ is stable. 

Let $v$ denote the vertex of an $m-$weighted tree $T$ and 
$e_1(v),\dots,e_n(v)$ the edges of $T$ starting at $v.$
The graph $T \setminus \{v,e_1(v),\dots,e_n(v)\}$ consists
of n weighted trees $T_1,\dots,T_n,$ which we call the subtrees
{\it complementary} to the vertex $v.$
A vertex $v$ of an m-weighted tree $T$  is called a {\it central vertex }
if 
\[wt(T_i) < \frac{m}{2}\] 

\noindent for every 
subtree complementary to $v.$
\bigskip

\begin{lema} \label{cent}
Let $T$ denote a stable $m-$weighted tree with the property \newline
(*) There is no edge $e$ of $T$ such that for the 2 complementary
subtrees $T_1$ and $T_2$ of $e$ we have $wt(T_1)=wt(T_2)=\frac{m}{2}.$ \newline
Then $T$ admits a unique central vertex $v_{cent}.$

\end{lema}

\bigskip
\noindent \textsc{Proof:I Existence of $v_{cent}.$} \newline
Start with any vertex $v_1$ of $T.$ Let $T_1^1,\dots,T_{n_1}^1$ denote
the subtrees complementary to $v_1.$
If $wt(T^1_j ) < \frac{m}{2}$ for all $j=1,\dots,n_1$ the vertex 
$v_1$ is central. Otherwise there is a unique subtree among the 
$T^{1}_{j}$ say $T_1^1$ such that $wt(T_1^1)=k > \frac{m}{2}.$
(Note that property (*) implies that $wt(T_1^1)$ cannot be equal
to $\frac{m}{2}).$ Let $v_2$ denote the other endpoint of the edge
$e_1(v).$
Let $T_1^2, \dots,T_{n_2}^2$ denote  the subtrees complementary to $v_2.$
We may assume that $T_1^2$ is the subtree
consisting of $v_1,T_2^1,\dots, T_{n_1}^1$
and the edges joining them.
Then
\[wt(T^2_1)=wt(v_1)+\sum_{\nu=2}^{n_1} wt(T_\nu^1)=m-k < \frac{m}{2}\]
On the other hand
\[ wt(v_2)+\sum_{\nu=2}^{n_2}wt(T^2_\nu)=wt(T_1^1)=k .\]
Since $wt(v_2) >0$ or $n_2 \geq 3$ we obtain that:
\[\max_{\nu=2}^{n_2} wt(T_\nu^2) <k.\]
If $\max_{\nu=2}^{n_2} wt(T_\nu^2)$ is still $> \frac{m}{2}$ we proceed
in the same way until we finally find a vertex
$v_{r}=v_{cent}$ with $wt(T^{r}_\nu) < \frac{m}{2}$ for all subtrees
$T_\nu^{r}$ complementary to $v_{r}.$

\noindent\textsc{Uniqueness of $v_{cent}.$}\newline
Suppose $v_{cent}^1 \neq v_{cent}^2$ are two central vertices.
$T$ being a tree, there is a unique path $v_1=v^1_{cent},v_2,v_3,
\dots,v_r=v^2_{cent}$ where $v_i$ and $v_{i+1}$ are connected
by one edge for $i=1,\dots,r-1.$
Let $T_1^1,\dots,T_{n_1}^1$ and $T_1^2,\dots,T_{n_2}^2$
denote the subtrees complementary to $v^1_{cent}$ and $v^2_{cent}.$
Without loss of generality we may assume
that $T_1^1$ is the subtree
consisting of $v_2,\dots,v_r,T_2^2,\dots,T_{n_2}^2$ and $T_1^2$ is the 
subtree consisting of $v_1,\dots,v_{r-1},T_2^1,\dots,T_{n_1}^1.$ 
Then:
\begin{eqnarray}
m & > & wt(T_1^1)+wt(T_1^2) \nonumber \\
  & \geq & \sum_{i=1}^{r}wt(v_i) + \sum_{\nu=2}^{n_1}wt(T_\nu^1)+
\sum_{\nu=2}^{n_2}wt(T^2_{\nu})                      \nonumber\\
  & = &wt(T)=m  \nonumber
\end{eqnarray}
and this is a contradiction. \qed

\noindent \textsc{Proof of Theorem ~\ref{main} :}
Let $C \in \overline\mathcal{M}_{0,m}$ be a stable $m-$marked 
curve of genus 0 for some $m \geq 3.$
Assume first that $C \not \in \Delta_{\frac{m}{2}}$ if $m$ is even.
If $T$ denotes the weighted dual graph associated to $C$ this means
that $T$ satisfies the condition (*) of Lemma ~\ref{cent}.
Hence, according to this Lemma, $C$ admits a central vertex
$v_{cent}$ dual to a {\it central component}
$C_{cent}.$ 
This means the following:
Let $C_1,\dots,C_r$ denote the weighted branches of the curve $C$
corresponding to the subtrees $T_1,\dots,T_r$
complementary to $v_{cent}.$
If $m_i$ denotes the number of marked points on $C_i$ then
$m_i < \frac{m}{2}.$
Choose coordinates $(x,y)$ of $C_{cent}\cong \P{1}$
and let $(x_1:y_1),\dots,(x_r:y_r)$ denote the points
of intersection
\[ (x_i:y_i)=C_i \cap C \mbox{ for } i=1,\dots ,r.\]
Moreover let $p_i=(x_i:y_i)$ for $i=r+1,\dots,s$ denote
the marked points of $C$ on the central component $C_{cent}.$
Note that the points $(x_i:y_i)$ are pairwise different
for $i=1,\dots,s.$ We then define:
\[ F_m(C):=\Pi_{i=1}^r (y_iX-x_iY)^{m_i} \Pi_{i=r+1}^s(y_iX-x_iY)\]
Geometrically the binary form $F_m(C)$ can be obtained
in the following way:
Contract the branches $C_i$ to the point $(x_i:y_i)$
for $i=1,\dots,r$ and associate to $(x_i:y_i)$
the weight  $wt(T_i),$ that is the number of marked points
of $C_i$ for $i=1,\dots,r.$
Then $F_m(C)$ is the binary form corresponding to the marked curve
$(C_{cent},(x_1:y_1)^{m_1},\dots,(x_r:y_r)^{m_r},(x_{r+1}:y_{r+1}),\dots,
(x_{s}:y_s)).$
This process can be extended to holomorphic families:
If $(\pi:\mathcal{C}\longrightarrow \mathcal{ U},\sigma_1,\dots,\sigma_m)$
denotes a holomorphic family of stable $m-$marked curves of genus
0 (that is $\sigma_1,\dots,\sigma_m$ are suitable sections of $\pi$)
it is easy to see that the central components of the fibres
form a holomorphic family $\mathcal{C}_{cent} \longrightarrow U$
and thus one obtains a holomorphic
family of binary forms 
$F_m (\mathcal{C}\stackrel{\pi}\longrightarrow U)$ over $\mathcal{U}.$
Thus we obtain a holomorphic map $F_m:\overline\mathcal{M}_{0,m} \longrightarrow
\overline\mathcal{B}_m$ for $m$ odd and $F_m:\overline\mathcal{M}_{0,m} 
\setminus \Delta_{\frac{m}{2}} \longrightarrow
\overline\mathcal{B}_m$ for $m$ even.
But then certainly $F_m$ extends to a continuous map on all
of $\overline\mathcal{M}_{0,m}$ by just mapping the divisor $\Delta_{\frac{m}{2}}$
to the semistable point of $\overline\mathcal{B}_m.$
Hence by Riemann's Extension Theorem $F_m$ is holomorphic everywhere.
\qed

\section{Local stable reduction } \label{sec4}

Consider the plane complete curve $C$ of arithmetic genus $g$ with affine 
equation:
\[y^2=(x-x_1)^{n_1}\dots(x-x_r)^{n_r} \]
with $x_i \neq x_j$ for $i \neq j.$
The curve $C$ admits only one point at $\infty$ namely $(0:1:0).$
If $n_0$ is its multiplicity, we have $\sum_{\nu=0}^{r} n_{\nu}=2g+2$
since $p_a(C)=g.$ Our aim in this section is to determine a stable 
model for $C.$

\bigskip

Let $C$ be a curve not necessarily reduced or irreducible. 
We assume that $C$ is the special fibre $\mathcal{C}_0$ of a fibration
$\pi:\mathcal{C} \longrightarrow U$ where $U$ 
is the unit disc and all fibers $\mathcal{C}_t=\pi^{-1}(t)$ are 
smooth for $t \neq 0$ and $\mathcal{C}_0=C.$

Recall (see \cite{kn:HM}) that a {\it reduction process of} $C$ (in 
$\mathcal{C}$)
consists of a finite sequence of steps of the form:\newline
\noindent i) A blow up of a point of $C.
$\\
\noindent ii) A base change $\mathcal{C}'$ by the $p^{th}$ power map
$U \longrightarrow U,z \mapsto z^p,$ where $p$ is a prime number, 
followed by the normalization  $n:\overline{\mathcal{C}} \longrightarrow 
\mathcal{C}':$ 
 \begin{center}
\begin{picture}(200, 100)
\put(1, 76){$\overline{\mathcal{C}} \stackrel{n}\longrightarrow$}
\put(42, 76){$\mathcal{C}'$}
\put(58,80){$\stackrel{{\phi}}{\vector(1, 0){95}}$} 
\put(30 , 55){$\pi'$}
\put( 45, 70){\vector(0,-1){35}}
\put(160,76){$\mathcal{C}$}
\put(162, 70){\vector(0, -1){35}}  
\put(168, 55){$\pi$}
\put(42, 17){$U$}
\put(58,20){$\stackrel{z \mapsto z^p}{\vector(1,0){95}}$}
\put(160,17){$U$}
\end{picture}
\end{center}

\noindent iii) Contraction of a smooth rational component
in the special fibre intersecting the other components
in at most 2 points.

The stable reduction theorem (see \cite{kn:HM}, p.118) says 
that for any curve $C$ there is a 
reduction process such that the resulting curve $\overline{C}$ is stable.
Notice that the map $\phi n: \overline\mathcal{C}:\longrightarrow \mathcal{C}$
is a covering map of degree $p$ branched exactly at those components
of the curve $C,$ whose multiplicity is not divisible by $p.$ 
In general there are many such coverings. 
We call step (ii) {\it the normalized base change of order p}.

If $p \in C$ is a point,
{\it local stable reduction 
of $C$} at $p$ is by definition a curve $\overline{C}$ such that:\newline
\noindent i)the curve $\overline{C}$ is obtained from $C$ by a  
{\it reduction process} with blow ups only at $p$ and its infinitely 
near points.\newline
\noindent ii) If $\overline{C}_n$ denotes the union of the new components,
that is those obtained by the sequence of blow ups, then:
\[ \overline{C}\setminus\overline{C}_n\cong C \setminus p \]
\noindent
iii) the curve $\overline {C}_n$ satisfies the stability condition
within $\overline{C}.$

Let now $C$ again denote the curve with affine equation as at the beginning 
of the section.
Choosing the coordinates appropriately we may assume: \newline
$n_0=
x_1=0.$
Hence C is of the form:

\begin{equation}\label{nf}
 y^2=x^{n_1}f(x)
\end{equation}

with $f(0) \neq 0$ and $n_1+deg(f)=2g+2.$

\begin{prop} \label{prop1} 
{\rm (a)}: A local  stable reduction of (\ref{nf}) with $n_1=2i$ is given by 
 $C' \cup E$ where $C'$ has affine equation $y^2=f(x)$
and $E$ has affine equation $y^2=z^{2i}-1.$ The curves $E$ and $C'$
intersect transversely in 2 points conjugate under the hyperelliptic
involution.\newline
{\rm (b)}: A local stable reduction of (\ref{nf}) with $n_1=2i+1$ is given by
$C' \cup E$ where $C'$ has affine equation $y^2=xf(x)$
and $E$ has affine equation $y^2=z^{2i+1}-1.$ The two components
intersect transversely in a point.
\end{prop}

Before we prove the proposition we need a lemma.

\begin{lema} \label{lema3}
Consider a part of a curve $C$ consisting  of three components $E_0,$
$E_1$ and $E_2$ such that $E_0$ intersects $E_1$ and $E_2$ 
transversely in one point 
and the multiplicities of the components $E_i$ are 
$n_0,n_1,$ and $n_2$ respectively.
Assume further that

\noindent i) The reduced curve $\overline{E}_i$ associated
to $E_i$ is of genus $0$ for $i=0,1,2.$\\
\noindent ii) $E_1$ and $E_2$ are the only components of $C$ intersecting $E_0.$\\
\noindent iii) $(n_0,n_i)=1,$ for $i=1,2.$

Then there is a reduction process of $C$ whose preimage of the above part of 
$C$ also consists of three curves $E_0,E_1$ and $E_2$ with $E_0$ intersecting
$E_1$ and $E_2$ transversely in one point 
but now the multiplicities of $E_0$ is $n_0=1$ while
the multiplicities $n_1$ of $E_1$ and $n_2$ of $E_2$ remain the same
and $E_0$ is of genus $0.$
 \end{lema}
\textsc{Proof:}
Let $p$ be a prime divisor of $n_0.$ According to ii)
$p \not \: \mid  n_i$ for $i=1,2.$ 
The $p^{th}$ power map (within a family of curves
$ \mathcal{C}\longrightarrow U$ as above)
 followed by normalization gives a $p:1$ covering
of $\overline{E}_0$ ramified at $\overline{E}_0 \cap \overline{E}_1$ and 
$\overline{E}_0 \cap \overline{E}_2.$
The multiplicity of the new $E_0$ is $n_0/p$ and Hurwitz formula implies
$g(E_0)=0.$ Repeating this process for the remaining prime divisors of $n_0/p$ 
yields the assertion. \qed

\bigskip

\textsc{Proof of Proposition \ref{prop1}:}(a): 
Blowing up i times the point $(0,0)$ gives the configuration of curves
indicated below. Here $E_j$ denotes the exceptional curve of the $(i+1-j)th$
blow-up. It is of multiplicity  $2j$ for $j=1,\dots,i$ indicated in the
picture as $\xymatrix@1{*+[F-]\txt{j}}$. 
The proper transform of $C$ is  denoted by $C'.$ It is given by the 
curve with affine equation $y^2=f(x).$
\vspace{.25in}

$\xymatrix{&&*+[F-]{2i}\ar@{-}[8,0]&&&&&&*+[F-]{4}\ar@{-}[4,0]&\\
            &&&&&&&*+[F-]{2}\ar@{-}[d];[3,0]\ar@{-}[0,2]&&E_1\\
             &*{C'\!\!\!\!\!\!\!\!\!\!\!\!\!\!}&&*{\!\!\!\!\!\!\!\!\!\!\
                 \!\!\!\!\!\!\!\!\!\!
                \!\!\!\
                \oval(120,70)[l]}&&&&*+[F-]{8}&&\\
&&&&&&*+[F-]{6}\ar@{-}[0,3]&&&E_3\\
          &&&&&&&E_4&E_2&\\
&&&&*+[F]{2i-4}\ar@{-}[3,0]&\ar@{.}[-1,1]&&&\\
&&&*+[F]{2i-6}\ar@{-}[0,6]&&&&&&E_{i-3}\\
&*+[F]{2i-2}  \ar@{-}[0,8]&&&&&&&&E_{i-1}\\
&&E_i&&E_{i-2}&&
}$

\vspace{.25in}

\noindent Making the normalized base change of order 2 we obtain 
the following configuration of curves with multiplicities as indicated.

\vspace{.25in}

$
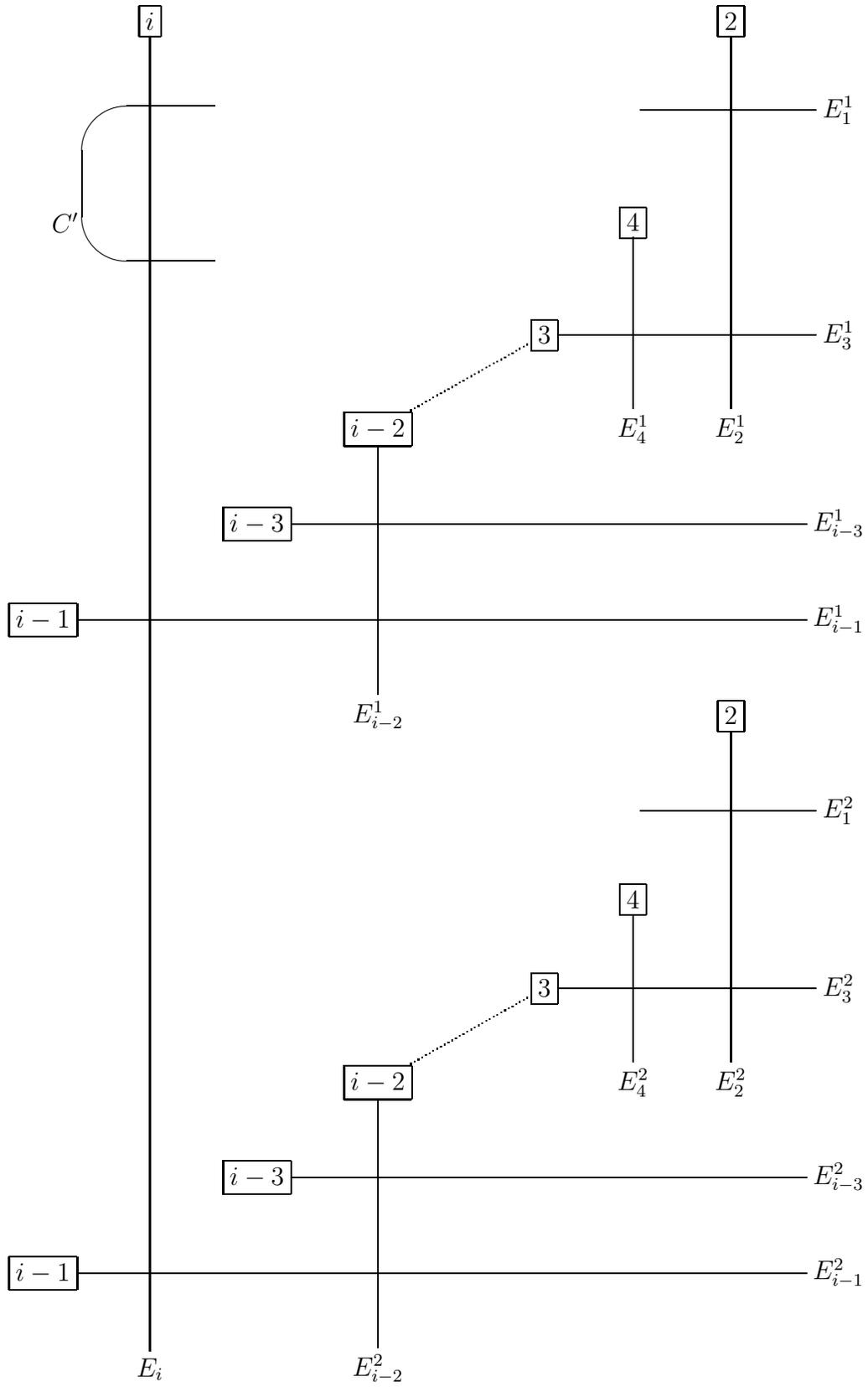
\xymatrix{&&*+[F-]{i}\ar@{-}[14,0]&&&&&&*+[F-]{2}\ar@{-}[4,0]&\\ \label{fig}
            &&&&&&&\ar@{-}[d];[3,0]\ar@{-}[0,2]&&E_1^1\\
             &*{C'\!\!\!\!\!\!\!\!\!\!}&&*{\!\!\!\!\!\!\!\!\!\!\
                 \!\!\!\!\!\!\!\!\!\!
                \!\!\!\
                \oval(120,70)[l]}&&&&*+[F-]{4}&&\\
&&&&&&*+[F-]{3}\ar@{-}[0,3]&&&E_3^1\\
&&&&*+[F]{i-2}\ar@{-}[3,0]\ar@{.}[-1,2]&&&E_4^1&E_2^1&\\
&&&*+[F]{i-3}\ar@{-}[0,6]&&&&&&E_{i-3}^1\\
&*+[F]{i-1}  \ar@{-}[0,8]&&&&&&&& E_{i-1}^1\\
&&&&E_{i-2}^1&&&&*+[F-]{2}\ar@{-}[4,0]&\\
            &&&&&&&\ar@{-}[d];[3,0]\ar@{-}[0,2]&&E_1^2\\
             &&&&&&&*+[F-]{4}&&\\
&&&&&&*+[F-]{3}\ar@{-}[0,3]&&&E_3^2\\
&&&&*+[F]{i-2}\ar@{-}[3,0]\ar@{.}[-1,2]&&&E_4^2&E_2^2&\\
&&&*+[F]{i-3}\ar@{-}[0,6]&&&&&&E_{i-3}^2\\
&*+[F]{i-1}  \ar@{-}[0,8]&&&&&&&& E_{i-1}^2\\
&&E_i&&E_{i-2}^2&&
}$
\vspace{.25in}

Consider the prime decomposition of i:

\[i=p_1\dots p_r\]
Then a)$p_j \not |\;\; (i-1),$ for $j=1,\dots,r$ and\\
b)If $p_j |k$ for some $k<i$ then $p_j \not | \;\; k \pm 1.$
 
Hence according to Lemma  \ref{lema3} the successive normalized base change
with $p_1,\dots,p_r$ yields a curve with the same configuration as above 
but where now the multiplicities are: 1 for $E_i,$ $i-1$ for $E_{i-1}^1,$
$i-1$ for $E_{i-1}^2$ and a divisor of $j$ for $E_j^k$ for $j=1,\dots,i-2,
k=1,2.$
In particular for all components $E_i^k, 2 \leq i \leq i-1,$ the conditions
of Lemma \ref{lema3} are satisfied. Applying normalized base change 
with the prime factors of the remaining multiplicities we obtain finally 
a curve with the same configuration where now all components,
which we denote by the same symbol, are reduced. 
All components 
$E_j^1$ and $E_j^2$ are rational, since they are cyclic coverings of rational 
curves ramified in 2 points. Hence we can successively contract 
$E_1^1,E_2^1,\dots,E_{i-1}^1$ and $E_1^2,\dots,E_{i-1}^2$ to obtain a curve
$\overline{C}=C' \cup E_i.$

The family of curves over $U$ obtained from
$\pi:\mathcal{C} \longrightarrow{U}$ by the composition 
of the base changes of orders $2,p_1,\dots,p_s$ is
a cyclic covering of degree $2i$ of $\mathcal{C}.$
Since its restriction to $E_i$ is unramified,
we obtain a cyclic covering $E_i \longrightarrow \p{1},$
where $\p{1}$ denotes the old component $E_i.$ 
Applying Hurwitz formula, we obtain $g(E_i)=g-1.$

Moreover we may assume that all fibres $\pi^{-1}(t),t \neq 0$
of the family $\pi: \mathcal{C} \longrightarrow U$ are hyperelliptic.
This implies that the curve $C' \cup E_i$ is hyperelliptic as a limit of 
hyperelliptic curves. Hence also the curve $E_i$ is hyperelliptic.
As  a cyclic covering of degree $2i$ of \p{1} it admits an automorphism of
 order 2i.
But it is well known (see e.g. \cite{kn:in}) that there is only one smooth
 hyperelliptic curve of genus $i-1$
with an automorphism of order $2i$ namely the curve with affine equation:
\[ y^2=z^{2i}-1.\] 
This concludes the proof.\qed

\bigskip

\textsc{Proof of Proposition \ref{prop1}}(b):
Blowing up $(i+2)$-times the point $(0,0)$ gives the configuration of curves:
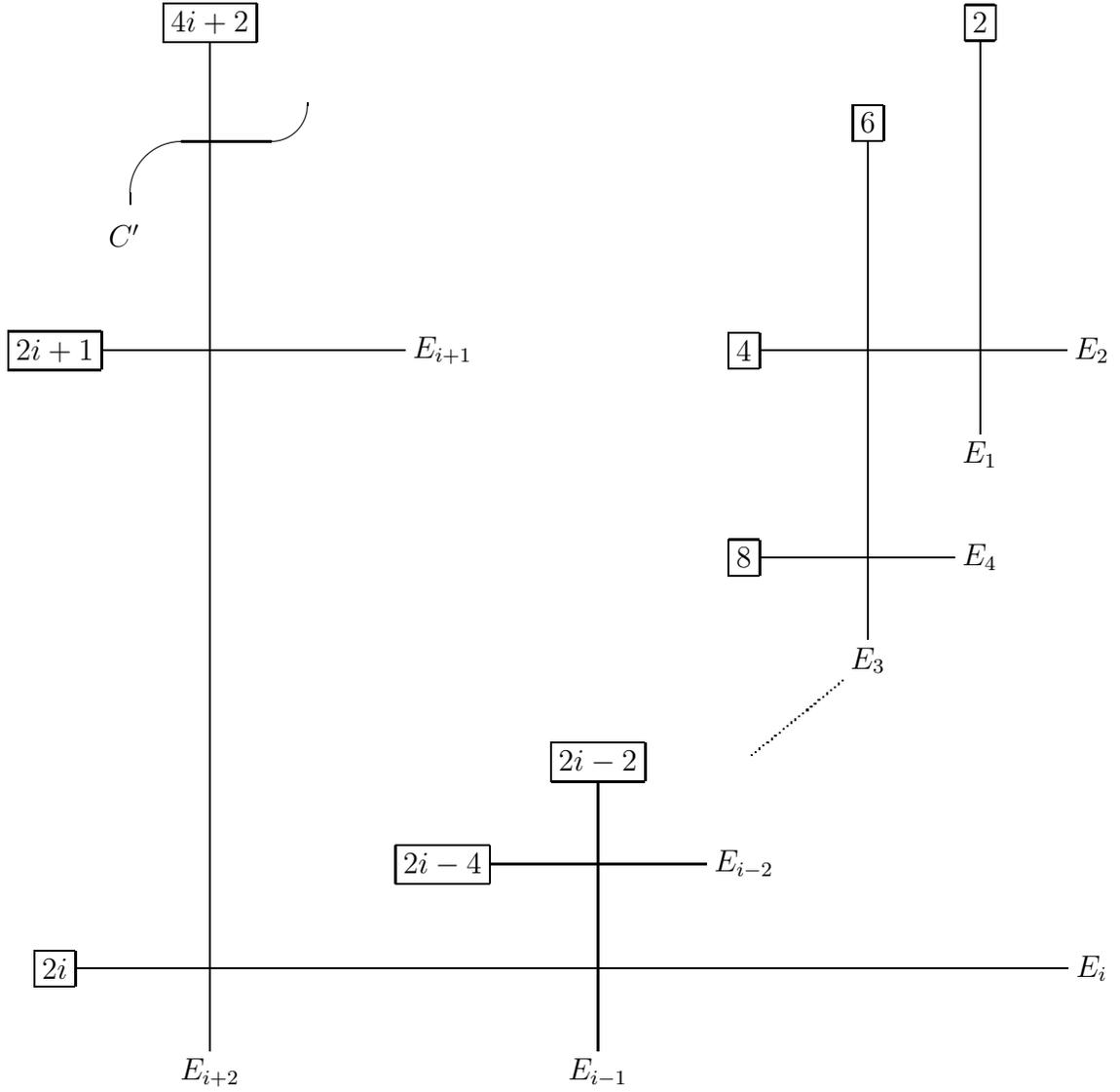
\vspace{.25in}

\[\xymatrix{&&*+[F-]{4i+2}\ar@{-}[10,0]&&&&&&*+[F-]{2}\ar@{-}[4,0]&\\
            &&*{\!\!\oval(80,30)[br]}&&&&&*+[F-]{6}\ar@{-}[5,0]&\\
             &*{C'\!\!\!\!\!\!\!\!\!\!\!\!\!\!\!\!\!\!\!\!\!\!\!\!\!\!\!}&&*{\!\!\!\!\!\!\!\!\!\!\!\!\!\!\!\!\!\
                                 \!\!\!\!\!\!\!\!\!\!\!
                \oval(100,48.5)[tl]}&&&&\\
&*+[F-]{2i+1}\ar@{-}[0,3]&&&E_{i+1}&&*+[F-]{4}\ar@{-}[0,3]&&&E_2\\
&&&&&&&&E_1&&\\
&&&&&&*+[F-]{8}\ar@{-}[0,2]&&E_4&&\\
          &&&&&&&E_3&&&&\\
&&&&&*+[F-]{2i-2}\ar@{-}[3,0]&\ar@{.}[-1,1]&&&\\
&&&&*+[F-]{2i-4}\ar@{-}[0,2]&&E_{i-2}&&&\\
&*+[F-]{2i}  \ar@{-}[0,8]&&&&&&&& E_{i}\\
&&E_{i+2}&&&E_{i-1}&&
}\]
\vspace{.25in}

\noindent Here $E_j$ denotes the exceptional curve of the (i+3-j)th blow up. 
It is of mutiplicity
$2j$ for $j=1,\dots,i$ of multiplicity $2i+1$ for $j=i+1$ and of multiplicity
 $4i+2$ for $j=i+2.$
The proper transform of $C$ is  denoted by $C'.$ It is given by affine 
equation $y^2=xf(x).$
Taking the normalized base change of order 2 we obtain the following 
configuration of curves with
multiplicities as indicated.

\vspace{.25in}

$
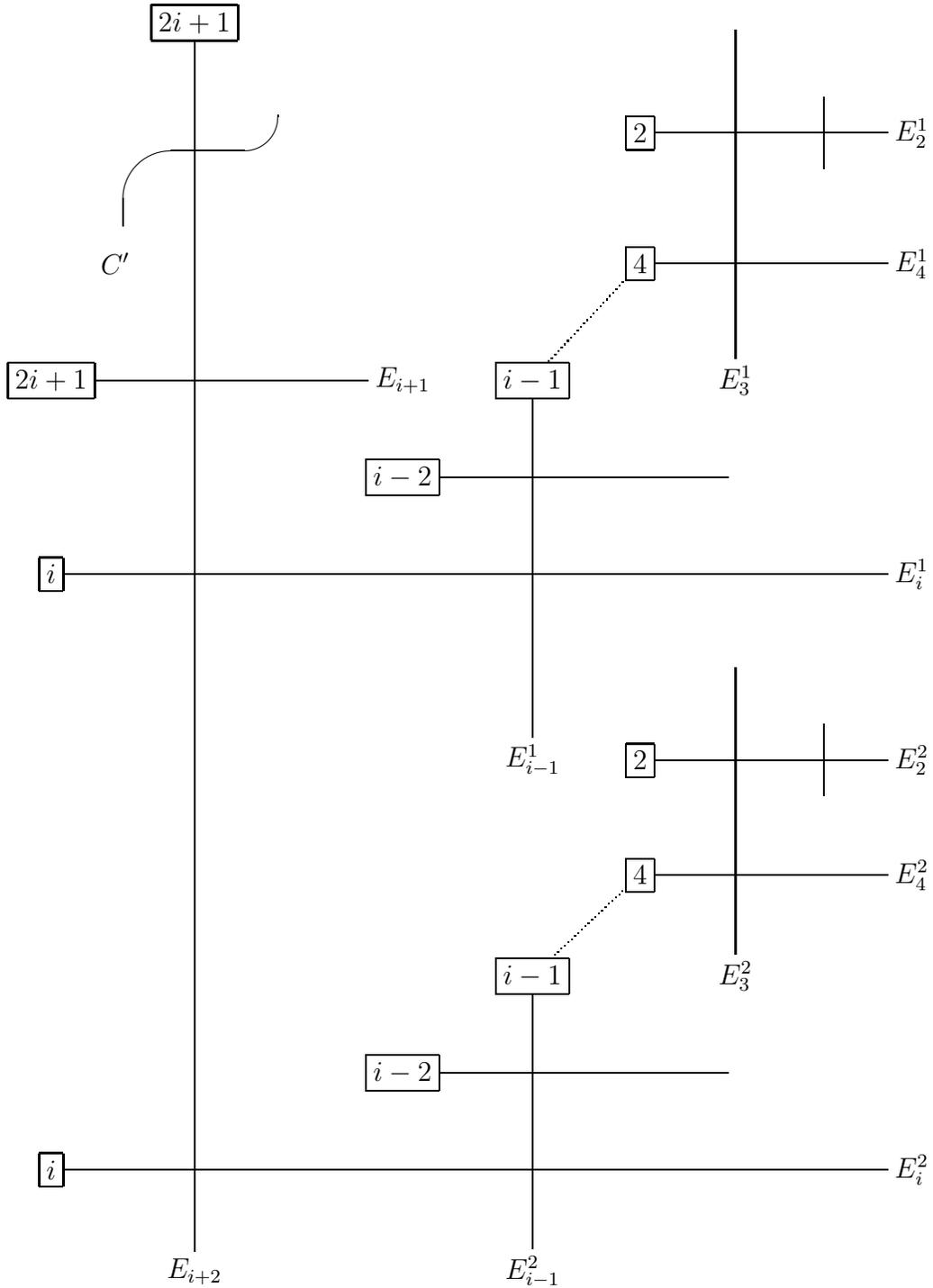
\xymatrix{&&*+[F-]{2i+1}\ar@{-}[12,0]&&&&&\ar@{-}[3,0]&&&\\
            &&*{\!\!\!\!\!\oval(80,30)[br]}&&&&*+[F-]{2}\ar@{-}[0,3]&
&*{\line(0,1){30}}&E_2^1\\
             &*{C'\!\!\!\!\!\!\!\!\!\!\!\!\!\!\!\!\!\!\!\!\!\!\!\!\!\!\!}&&*{\!\!\!\!\!\!\!\!\!\!\!\!\!\!\!\!\!\
                                 \!\!\!\!\!\!\!\!\!\!
                \oval(100,63.2)[tl]}&&&*+[F-]{4}\ar@{-}[0,3]&&&E_4^1\\
&*+[F-]{2i+1}\ar@{-}[0,3]&&&E_{i+1}&*+[F-]{i-1}\ar@{-}[4,0]\ar@{.}[-1,1]&
                                                                    &E_3^1&&\\
&&&&*+[F-]{i-2}\ar@{-}[0,3]&&&&&&\\
&*+[F-]{i}  \ar@{-}[0,8]&&&&&&&& E_{i}^1\\
&&&&&&&\ar@{-}[3,0]&&&&\\
            &&&&&E_{i-1}^1&*+[F-]{2}\ar@{-}[0,3]&&*{\line(0,1){30}}&E_2^2\\
             &&&&&&*+[F-]{4}\ar@{-}[0,3]&&&E_4^2\\
&&&&&*+[F-]{i-1}\ar@{-}[3,0]\ar@{.}[-1,1]&&E_3^2&&\\
&&&&*+[F-]{i-2}\ar@{-}[0,3]&&&&&&\\
&*+[F-]{i}  \ar@{-}[0,8]&&&&&&&& E_{i}^2\\
&&E_{i+2}&&&E_{i-1}^2&&
}$

\vspace{.25in}

Consider the prime decomposition of $2i+1.$

\[2i+1=p_1\dots p_r\]
Then\\
a) $p_j  \not |\;\; i$ for $j=1,\dots r.$ \\
b) If $p_j\;\; |\;\; k$ for some $k<i,$ then $p_j \not \: \mid \;\; k \pm 1.$

This implies that according to Lemma \ref{lema3} the successive 
normalized base changes with order $p_1,\dots p_r$ yield a curve with the 
same configuration apart from the fact that the curve $E_{i+1}$
is replaced by $2i+1$ curves $E_{i+1}^1,\dots, E_{i+1}^{2i+1}$ intersecting
$E_{i+2}$ only once and not intersecting any other component.
The curves $C',E_{i+1}^1,\dots, E_{i+1}^{2i+1}$ are of multiplicity 1 
and the multiplicity of $E_{\nu}^1$ and
$E_{\nu}^2$ is a divisor of $\nu$ for $\nu=1,\dots,i.$
Hence to the remaining multiplicities one can also apply Lemma \ref{lema3} to obtain a curve
with the same configuration but where now all components are reduced.

According to Lemma \ref{lema3} all curves $E_{\nu}^1,E_{\nu}^2$ for 
$\nu=1,\dots,i$ are rational.
Moreover $E_{i+1}^1,\dots E_{i+1}^{2i+1}$ are rational. Hence we can 
successively contract $E_1^1,E_2^1,\dots,E_i^1,E_1^2,\dots,E_i^2$ 
and $E_{i+1}^1,\dots,E_{i+1}^{2i+1}$ to obtain the stable curve 
$\overline{C}=C' \cup E_{i+2}.$
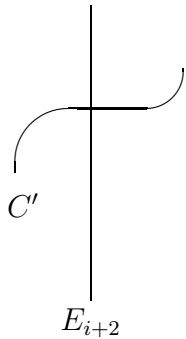
\vspace{.25in}   
\[\xymatrix{&&\ar@{-}[3,0]&&&&&&&&\\
            &&*{\!\!\!\!\!\oval(80,30)[br]}&&&&&&\\
             &*{C'\!\!\!\!\!\!\!\!\!\!\!\!\!\!\!\!}&&*{\!\!\!\!\!\!\!\!\!\!\!\!\!\!\!\!\!\
                                 \!\!\!\!\!
                \oval(100,48.5)[tl]}&&&&&&\\
&&E_{i+2}&&&&&&&
}\]
\vspace{.25in} 

Hurwitz formula yields $g(E_{i+2})=i.$ Moreover in the same way as in the
 proof of Proposition~\ref{prop1}, 
one can see that $E_{i+2}$ is a hyperelliptic curve of genus $i$ admitting 
an automorphism 
of order $2i+1.$ But it is well known (see e.g. \cite{kn:in}) that there is 
only one such curve 
namely the curve with affine equation
\[ y^2=z^{2i+1}-1,\]
thus terminating the proof.\qed

\section{The holomorphic map $f_g: \g \longrightarrow \binf$}

In Corollary~\ref{iso1}, we saw that there is a canonical isomorphism
 $ \Phi :\overline{\mathcal{H}}_g \longrightarrow \overline{\mathcal{M}}_{0,2g+2}$
and Theorem~\ref{main} says that the canonical map 
$\mathcal{M}_{0,2g+2}\longrightarrow \bin$ extends to a holomorphic map 
$F_{2g+2}: \overline{\mathcal{M}}_{0,2g+2}\longrightarrow \binf.$
Since the composition $f_g:=F_{2g+2}\circ \Phi$ certainly extends 
the canonical isomorphism $\mathcal{H}_g \widetilde{\longrightarrow} \bin,$
this completes the proof of Theorem~\ref{princ}.

\bigskip

Using the results of Section~\ref{sec4}, the map $f_g$ can also be described
as follows: let $C$ be a stable hyperelliptic curve and 
$C'\longrightarrow B$ the associated admissible double cover
(see Proposition~\ref{iso}). If $C \not \in \Delta_g,$ then according to
Theorem~\ref{main} the curve $B$ admits a unique central component
$B_{cent}.$ Let $B_1, \dots, B_r$ denote the closure in $B$ of the connected
components of $B \setminus B_{cent}$ and $x_i:=B_i \cap B_{cent}$ 
for $i=1,\dots,r.$ 
 The components $B_i$ are themselves
models of hyperelliptic curves, say  of genus $g_i$ and this can be deformed
to the curves $E_i$ of Proposition~\ref{prop1}. According to this proposition
the curves $E_i$ can be contracted. If $x_{r+1},\dots,x_s$ are the smooth
ramification points of $B_{cent}=\p{1},$
the map $f_g$ associates to $C$ the binary form:
\[ f_g(C)=(X-x_1Y)^{n_1}\dots(X-x_sY)^{n_s}\]

with
\[ n_i= \left\{ \begin{array}{cccccc}
         2g_i+2&&&&i\leq r&\mbox{ and $g_i$ even }\\
         2g_i+1&&\mbox{if}&&i \leq r&\mbox{ and $g_i $ odd}\\
            1&&&&i>r&
\end{array}\right. \]

Finally, we want to study the behaviour of the map $f_g$
on the boundary divisors $\Delta_i$ and $\Xi_i.$
\begin{prop}
The holomorphic map $f_g$ contracts $\Delta_i$ (respectively
$\Xi_i$) to a subvariety of dimension $2g-2i-1$ (respectively
$2g-2i-2),$ except when  $2i+1=g+1$ (resp. $2i+2=g+1$). In this case, 
 $f_g$
contracts $\Delta_i$ (resp. $\Xi_i$) to a point.
\end{prop}

\textsc{Proof:}
Let $C$ be a general element of $\Delta_i$ and 
$C \longrightarrow B$ the associated double cover.
Then $B=B_1 \cup B_2$
where $B_1,B_2$ are irreducible components of genus zero.
If $2i+1=g+1=2g+2/2$ then $f_g$ maps $C$ to the semistable point.
Otherwise and
without loss of generality we can assume that $B_1$ is the central
component and therefore it has $2g-2i+1$ marked points,
say $x_1,\dots,x_{2g-2i+1}.$
If $x_0$ denotes the point of intersection
of $B_1$ and $B_2,$ then according to the above description 
of the map $f_g$ the image $f_g(C)$ is the binary form:
\[f_g(C)=(X-x_0)^{2i+1}(X-x_1)\dots(X-x_{2g-2i+1}).\]
It is clear that every binary form of this type is contained
in the image $f_g(\Delta_i).$
Moreover, the forms of this type make up an open set
set of the variety $f_g(\Delta_i).$
This implies 
\[ dim(f_g(\Delta_i))=2g-2i+2-dim(Aut(\p{1}))=2g-2i-1,\]
since $dim(Aut(\p{1}))=PGL_1(\complex{})$ is of dimension 3.
The computation of $dim(\Xi_i)$ is analogous. \qed

\begin{center}
\vskip20pt
$\begin{array} {c}
\hbox{\footnotesize Departamento de Matem\'atica, UFMG}\\
\hbox{\footnotesize Belo Horizonte, MG 30161--970, Brasil.}
\\
\hbox{\footnotesize dan@mat.ufmg.br}
\end{array}$\quad\quad
$\begin{array} c

\hbox{\footnotesize Mathematisches Institut}\\
\hbox{\footnotesize Bismarckstr. $1\frac{1}{2},$ 91054, Erlangen}
\\

\hbox{\footnotesize lange@mi.uni-erlangen.de}
\end{array}$
\end{center}
\end{document}